\documentclass{amsart}
\usepackage{color,graphicx,wrapfig}
\usepackage[notref,notcite]{showkeys}
\newcommand{\non}{\nonumber}
\usepackage{graphicx}
\usepackage{amsmath}
\usepackage{amsfonts}
\usepackage{amssymb}
\usepackage{amsthm}
\usepackage{txfonts}
    \def\haus{\,d{\mathcal H}^{n-1}}
    \def\dh1{{\,d{\mathcal H}^1}}
    \def\dh2{{\,d{\mathcal H}^2}}
    
    \def\R{{\mathbb{R}}}\def\N{{\mathbb{N}}}
    \def\Hf{{\mathbb{H}}}

    \def\e{\mathbf{e}}
    \def\dist{\mathbf{dist}}
    \def\d{{\boldsymbol{\delta}}}

    \def\u{{\mathbf u}}
    \def\v{{\mathbf v}}

\def\z{{\mathbf z}}
\def\d{{\mathbf d}}
\def\f{{\mathbf f}}
\def\vphi{{\boldsymbol\phi}}
\def\c{{\mathbf c}}
\def\h{{\mathbf h}}
\def\g{{\mathbf g}}
\def\w{{\mathbf w}}
\def\p{{\mathbf p}}
\def\q{{\mathbf q}}
\def\a{{\mathbf a}}

 \def\supp{\mathbf{supp}\textrm{ }}

\newtheorem{lem}{Lemma}
\newtheorem{prop}{Proposition}
\newtheorem{cor}{Corollary}
\newtheorem{thm}{Theorem}
\newtheorem{definition}{Definition}

\newtheorem{remark}{Remark}

\title{Equilibrium points of a singular cooperative system with free boundary}
\author[J. Andersson]{John Andersson}
\address{Department of Mathematics, Royal Institute of Technology,
100~44  Stockholm, Sweden}
\email{johnan@kth.se}
\author[H. Shahgholian ]{Henrik Shahgholian}
\address{Department of Mathematics, Royal Institute of Technology,
100~44  Stockholm, Sweden}
\email{henriksh@math.kth.se}
\urladdr{http://www.math.kth.se/~henriksh/}
\author[N. Uraltseva]{Nina N. Uraltseva}
\address{St-Petersburg State University,
Universitetsky pr. 28, Stary Petergof, 198504, RUSSIA}
\email{uraltsev@pdmi.ras.ru}
\author[G.S. Weiss]{Georg S. Weiss}
\address{Department of Mathematics, Heinrich Heine University, 40225 D\"usseldorf}
\email{weiss@math.uni-duesseldorf.de}
\thanks{$2000$ {\it Mathematics Subject Classification.\/} Primary
35R35, Secondary  35J60}
\thanks{H. Shahgholian has been supported in part by
the Swedish Research Council. Nina Uraltseva was supported by Russian Foundation of Basic research (RFBR) Grant 11-01-00825, and by Grant of St-Petersburg State University 6.38.670.2013.
Both G.S. Weiss and N. Uraltseva thank the G\"oran Gustafsson Foundation
for visiting appointments to KTH}
\date{\today}

\begin{document}

\maketitle

\begin{abstract} In this paper we initiate the study of maps minimising the energy 
  $$
    \int_{D} (|\nabla \u|^2+2|\u|)\  dx.
  $$
which, due to Lipschitz character of the integrand, gives rise to the singular Euler equations
  $$
  \Delta \u=\frac{\u}{|\u|}\chi_{\left\lbrace |\u|>0\right\rbrace }, \qquad \u = (u_1, \cdots , u_m) \ .
  $$
Our primary goal in this paper is to set up a road map for future developments of the theory related to such energy minimising maps. 

Our results here concern regularity of the solution as well as that of the free boundary. They are achieved by using 
 monotonicity formulas and epiperimetric inequalities, in combination with geometric analysis.  
\end{abstract}


\section{Introduction}


\subsection{Background}

In this paper we shall study the singular system
\begin{equation}\label{system}
\begin{array}{l}
\Delta \u=\frac{\u}{|\u|}\chi_{\left\lbrace |\u|>0\right\rbrace }, \qquad \u = (u_1, \cdots , u_m) \ ,
\end{array}
\end{equation}
where $\u: \R^n \supset D \to \R^m$, $n\ge 2$, $m\geq 1$, and $|\cdot|$ is the Euclidean norm on the respective
spaces.
System (\ref{system}) is a particular example of the equilibrium state of a cooperative system:
the corresponding reaction-diffusion system
\begin{align*}
&u_t - \Delta u = -\frac{u}{\sqrt{u^2+v^2}},\\
&v_t - \Delta v = -\frac{v}{\sqrt{u^2+v^2}}
\end{align*}
would mean that, considering the concentrations $u$ and $v$ of
two species/reactants, each species/reactant
slows down the extinction/reaction of the other species.
The special choice of our reaction kinetics would assure
a {\em constant} decay/reaction rate in the case that $u$ and $v$ are
of comparable size.  

System (\ref{system}) may also be seen as one of the simplest extensions
of the classical {\em obstacle problem} to the vector-valued case:
Solutions of the classical obstacle problem are minimisers of the energy
$\int_{D} (\frac{1}{2}|\nabla u|^2 + \max(u,0))\  dx,$
where $u: \R^n \supset D \to \R$.
Solutions of (\ref{system}) are minimisers of the energy
\begin{equation}\label{main-functional}
  \int_{D} (|\nabla \u|^2+2|\u|)\  dx.
\end{equation}

It is noteworthy that in the scalar case, i.e. when $m=1$, one recovers the two phase free boundary problem 
$$
\Delta u = \chi_{\{ u>0\}} -  \chi_{\{u<0\}}
$$
contained in the analysis of \cite{monats}. While \cite{monats} as well as
the two-phase result \cite{cks} relied essentially
on the use of the monotonicity formula by Alt-Caffarelli-Friedman \cite{acf},
a corresponding formula seems to be unavailable in our vector-valued
problem.

There are several results concerning the obstacle problem for systems of various types: Optimal switching, 
multi-membranes, control of systems, constrained  weakly elliptic systems, vector-valued obstacle problems, and probably many others.
Although  not directly relevant to our work, we refer to some papers that might be of interest for the readers 
\cite{ad-i}, \cite{ad-ii}, \cite{cv}, \cite{ef}, \cite{fuchs}, \cite{duz-fuchs}, \cite{lions}.

\subsection{Main Result and Plan of the paper} In this paper we are interested in qualitative behavior of the minimisers $\u$ of the functional  (\ref{main-functional}) as well as of the free boundary $\partial \{ x: |\u(x)|>0\}$;
here  $\u = (u_1, \cdots, u_m)$ and $m \geq  1$.
Note that the part of the free boundary where the gradient $\nabla \u \ne 0$,
is by the implicit function theorem locally a $C^{1,\beta}$-surface, so
that we are more concerned with the part where the gradient vanishes.

The main results of this paper  (presented in Theorem \ref{regular})
states that the set of "regular" free boundary points of the  minimisers $\u$ to the functional  (\ref{main-functional}) are locally a $C^{1,\beta}$
surface.

In proving this result we need an array of technical tools including monotonicity formulas (Lemma \ref{mon} in Section \ref{mon-sec}), 
quadratic growth of solutions (Theorem \ref{quadratic}), and an epiperimetric inequality  (Theorem \ref{epi}), for the balanced energy functional (\ref{M-functional}).

An epiperimetric inequality has been proved in \cite{inv} 
by one of the authors for the scalar obstacle problem. See also \cite{monneau} for a 
related approach to the scalar obstacle problem with Dini
continuous coefficients.



\subsection{Notation}
Throughout this paper $\R^n,\R^m,\R^{nm}$ etc. will be equipped with the Euclidean
inner product $x\cdot y$ and the induced norm $\vert x \vert\> ,\>
B_r(x^0)$ will denote the open $n$-dimensional ball of center
$x^0\> ,$ radius $r$ and volume $r^n\> \omega_n\> ,
\> B'_r(x^0):= \{x\in B_r(x^0)\> : \>
x_n = (x^0)_n\}\> , \> B^+_r(x^0) := \{x\in B_r(x^0)\> : \>
x_n > (x^0)_n\}$ and $\e^i$ the $i$-th unit vector in $\R^k$.
If the center $x^0$ is not specified, then it is assumed to be the origin.
Given a set $A\subset \R^n\> ,$ we denote its interior by $A^\circ$
and its characteristic function by $\chi_A\> .$
In the text we use the $n$-dimensional Lebesgue-measure
$|A|$ of a set $A$ and the $k$-dimensional Hausdorff-measure
$\mathcal{H}^k\> .$ When considering the boundary of a given set,
$\nu$ will typically denote the topological outward
normal to the boundary and $\nabla_\theta f := \nabla f
\> - \> \nabla f \cdot \nu \> \nu$ the surface derivative
of a given function $f\> .$ Finally, we shall often use
abbreviations for inverse images like $\{u>0\} := 
\{x\in D\> : \> u(x)>0\}\> , \> \{x_n>0\} := 
\{x \in \R^n \> : \> x_n > 0\}$ etc. and occasionally 
we employ the decomposition $x=(x',x_n)$ of a vector $x\in \R^n\> .$
Last, let $\Gamma(\u):= D\cap \partial \{x\in D: |\u (x)| >0\}$ and $\Gamma_0(\u):= \Gamma(\u) \cap  \{x: \nabla \u (x)= 0\} $.


\section{The Epiperimetric Inequality}


Following \cite{inv}, we prove in this section  an epiperimetric inequality, which 
tells us that close to half-plane solutions, the minimal energy achieved is
lower than that of $2$-homogeneous functions, and the energy difference can be estimated.
This will imply in later sections a certain non-degeneracy of
the energy close to half-plane solutions, and ultimately lead to regularity of the free boundary.
Since the epiperimetric inequality is rather an abstract property of the energy,
and represents the core of our result,
we put this section at the beginning. Although the proof follows partly
the proof in \cite{inv}, the PDE resulting from the ``linearization'' carried
out in the proof is different from that in \cite{inv} and introduces new difficulties.

Let 
\begin{equation}\label{M-functional}
  M(\v):= \int_{B_1} (|\nabla \v|^2 + 2 |\v|) -  2 \int_{\partial B_1} |\v|^2 \haus ,
\end{equation}
and let 
$$\Hf := \{ \frac{\max(x\cdot \nu,0)^2}{2} \e : 
\nu \textrm{ is a unit vector in }\R^n \textrm{ and }
\e \textrm{ is a unit vector in }\R^m\}.$$
We define $$\frac{\alpha_n}{2}:= M(\frac{\max(x\cdot \nu,0)^2}{2} \e).$$
\begin{thm}\label{epi}
There exists $\kappa\in (0,1)$ and $\delta>0$ such that if
$\c$ is a homogeneous
function of degree $2$ satisfying $\Vert \c-\h\Vert_{W^{1,2}(B_1;\R^m)}+\Vert \c-\h\Vert_{L^\infty(B_1;\R^m)}\le \delta$ for some $\h \in \Hf$,
then there is a $\v\in W^{1,2}(B_1;\R^m)$ such that $\v=\c$ on $\partial B_1$ and
\begin{equation}\label{epiin}
M(\v) \le (1-\kappa) M(\c) + \kappa \frac{\alpha_n}{2}.
\end{equation}
\end{thm}
{\it Remark:} Note that the closeness in $L^\infty$ is not really necessary and is assumed only 
in order to avoid capacity arguments in the proof.

{\it Proof of the Theorem.}
Suppose towards a contradiction that there are sequences $\kappa_k\to 0,
\delta_k\to 0$, $\c_k\in W^{1,2}(B_1; \mathbb{R}^m)$,  and $\h_k \in \Hf$ 
such that $\c_k$ is homogeneous function of degree $2$ and satisfies 
$$
\Vert \c_k-\h_k\Vert_{W^{1,2}(B_1;\R^m)}= \delta_k, \qquad \Vert \c_k-\h_k\Vert_{L^\infty(B_1;\R^m)}\to_{k\to \infty } 0
$$
and that
\begin{equation}\label{epi1}
M(\v) > (1-\kappa_k) M(\c_k) + \kappa_k \frac{\alpha_n}{2}
\qquad \textrm{ for all }   \v\in \c_k +  W_0^{1,2}(B_1;\R^m) .
\end{equation}
Rotating in $\R^n$ and in $\R^m$ if necessary we may assume that 
$$\h_k(x) = \frac{\max(x_n,0)^2}{2} \e^1 =:\h$$
where $\e^1=(1,0,\dots,0)\in \R^m$. Subtracting from (\ref{epi1})
$M(\h)=\frac{\alpha_n}{2}$, we obtain
\begin{equation}\label{epi2}
(1-\kappa_k) (M(\c_k) -M(\h)) < M(\v)-M(\h)
\textrm{ for every } \v\in W^{1,2}(B_1;\R^m) \textrm{ such that }
\v=\c_k \textrm{ on }\partial B_1.
\end{equation}
Observe now that for each $\phi\in W^{1,2}(B_1)$ and $h := \frac{\max(x_n,0)^2}{2}$
$$
2\int_{B_1} (\nabla h\cdot \nabla \phi+ \chi_{\{x_n>0\}}\phi)
- 2 \int_{\partial B_1} 2h\phi\haus = 
2 \int_{\partial B_1} (\nabla h\cdot \nu - 2h)\phi\haus
= 0,$$
and therefore
$$I := 2 \int_{B_1} (\nabla \h\cdot \nabla (\c_k-\h)+ \chi_{\{x_n>0\}} \e^1\cdot
(\c_k-\h))
- 2 \int_{\partial B_1} 2\h\cdot  (\c_k-\h)\haus=0.$$
Subtracting $(1-\kappa_k)I$ from the left-hand side of (\ref{epi2})
and subtracting $I$ with $\c_k$ replaced by $\v$ from the right-hand side of (\ref{epi2}),
we obtain thus
\begin{align*}
&(1-\kappa_k) \bigg[\int_{B_1} (|\nabla \c_k|^2 + 2 |\c_k|) - 2 \int_{\partial B_1} 
|\c_k|^2 \haus -\int_{B_1} (|\nabla \h|^2 + 2 |\h|) - 2 \int_{\partial B_1} 
|\h|^2 \haus\\&
- 2 \int_{B_1} (\nabla \h\cdot \nabla (\c_k-\h)+ 2 \chi_{\{x_n>0\}} \e^1\cdot
(\c_k-\h))
+ 2 \int_{\partial B_1} 2\h\cdot  (\c_k-\h)\haus\bigg]
\\&<
\int_{B_1} (|\nabla \v|^2 + 2 |\v|) - 2 \int_{\partial B_1} 
|\v|^2 \haus -\int_{B_1} (|\nabla \h|^2 + 2 |\h|) - 2 \int_{\partial B_1} 
|\h|^2 \haus\\&
- 2 \int_{B_1} (\nabla \h\cdot \nabla (\v-\h)+ 2 \chi_{\{x_n>0\}} \e^1\cdot
(\v-\h))
+ 2 \int_{\partial B_1} 2\h\cdot  (\v-\h)\haus.
\end{align*}
Rearranging terms yields
\begin{align}\label{epi3}
&(1-\kappa_k) \bigg[\int_{B_1} |\nabla (\c_k-\h)|^2  - 2 \int_{\partial B_1} 
|\c_k-\h|^2 \haus
+ 2 \int_{B_1^-} |\c_k| + 2  \int_{B_1^+} (|\c_k|- \e^1\cdot \c_k)\bigg]
\\&<
\int_{B_1} |\nabla (\v-\h)|^2 - 2 \int_{\partial B_1} 
|\v-\h|^2 \haus
+ 2 \int_{B_1^-} |\v| + 2  \int_{B_1^+} (|\v|- \e^1\cdot \v).\non
\end{align}
Define now the sequence of functions  $\w_k := (\c_k-\h)/\delta_k$. 
Then  $\Vert \w_k\Vert_{W^{1,2}(B_1;\R^m)}=1$ 
and, passing to a subsequence if necessary,
$\w_k \to \w$ weakly in $W^{1,2}(B_1;\R^m)$.
In order to obtain a contradiction, we are going to prove that
$\w_k \to \w$ strongly in $W^{1,2}(B_1;\R^m)$ and that $\w\equiv 0$ in $B_1(0)$.
\\
{\bf Step 1:
$\w\equiv 0$ in $B_1^-$, and
$\int_{B_1^+} (|\c_k|- \e^1\cdot \c_k) \le C \delta_k^2$}

Plug in $\v := (1-\zeta) \c_k + \zeta \h$ in (\ref{epi3}),
where $\zeta\in W^{1,2}_0(B_1)$ is radial symmetric and satisfies $0\le\zeta\le 1$.
Since $(\v-\h)/\delta_k = (1-\zeta) \w_k$,
we obtain
\begin{align*}
&2(1-\kappa_k) 
\left[\int_{B_1^-} \frac{|\c_k|}{\delta_k^2} + \int_{B_1^+} \frac{|\c_k|- \e^1\cdot \c_k}{\delta_k^2}\right]
\\&< C_1 
+ 2 \int_{B_1^-} (1-\zeta)\frac{|\c_k|}{\delta_k^2} + 2  \int_{B_1^+} (1-\zeta)
\frac{|\c_k|- \e^1\cdot \c_k}{\delta_k^2}
\end{align*}
and
\begin{align*}&
\int_{B_1^-} (\zeta-\kappa_k)\frac{|\c_k|}{\delta_k^2} + \int_{B_1^+} 
(\zeta-\kappa_k) \frac{|\c_k|- \e^1\cdot \c_k}{\delta_k^2}
\le C_1.\end{align*}
Using the homogeneity of $\c_k$ we see that for large $k$,
\begin{align*}
&\int_{B_1^-} (\zeta-\kappa_k)|\c_k|
= \int_0^1 (\zeta(\rho)-\kappa_k)\rho^{n+1}\> d\rho  \int_{\{ x_n<0\} \cap \partial B_1}|\c_k|\haus 
\ge c_0 \int_{\{ x_n<0\} \cap \partial B_1}|\c_k|\haus,\end{align*}
where $c_0>0$ depends only on $\zeta$ and $n$.
We also get the corresponding estimate in $B_1^+$.
It follows that
\begin{equation}\label{epi4}
\int_{B_1^-} |\c_k|\le C_2 \delta_k^2 \textrm{ and that }
\int_{B_1^+} \left( |\c_k|- \e^1\cdot \c_k \right) \le C_2 \delta_k^2. 
\end{equation}
In particular, $$\int_{B_1^-} |\w_k|\le C_2 \delta_k,$$
implying the statement of Step 1.
\\
{\bf Step 2: } $\Delta(\e^1\cdot \w )=0$ in $B_1^+(0), \e^j\cdot\w=d_jh$ in
$B_1^+(0)$  for each $j>1$, and some constant  $d_j$.

Fix a  ball $B\subset\subset B_1^+$ and plug 
$\v := (1-\zeta) \c_k + \zeta (\h+\delta_k \g)$ into (\ref{epi3}),
where $\zeta\in C^\infty_0(B_1^+)$ and $\g\in W^{1,2}(B_1;\R^m)$
such that $\zeta\equiv 1$ in $B$, $\zeta\equiv 0$ in $B_1^-$
and $\g $ is a bounded $W^{1,2}(B_1;\R^m)$-function.
Observing that
$$\frac{\v-\h}{\delta_k} = (1-\zeta)\frac{\c_k-\h}{\delta_k} + \zeta\g,$$
we obtain ---using (\ref{epi4}) as well as the fact that $\supp \zeta \subset\subset B_1^+$--- that
\begin{align*}
&\int_{B_1^+}(2\zeta - \zeta^2) |\nabla \w_k|^2 + \frac{2}{\delta_k^2} \int_{B_1^+}
 \zeta(|\c_k|-\e^1\cdot \c_k)\\ &
\le o(1)+\int_{B_1^+}\zeta^2|\nabla \g|^2 + 2 \int_{B_1^+\setminus B} 
\left (|\nabla \zeta|^2|\g-\w_k|^2+(2\zeta-2\zeta^2)\nabla \w_k\cdot \nabla \g\right) \\
&+2\int_{B_1^+\setminus B} \left((1-\zeta)\nabla \zeta\cdot \nabla \w_k(\g-\w_k)+
\zeta\nabla \zeta \cdot \nabla \g(\g-\w_k)\right)\\
& +\frac{2}{\delta_k^2} \int_{B_1^+}
\zeta\left(|\h+\delta_k\g|-\e^1\cdot (\h+\delta_k \g)\right).
\end{align*}

Note that $\delta_k \w_k \to 0$ uniformly in $B_1$. Therefore we have on $\supp \zeta$
 $$|\c_k|-\e^1\cdot \c_k=
(h+\delta_k \e^1\cdot \w_k) \left(\sqrt{1+  \delta_k^2\frac{|\w_k|^2-(\e^1\cdot \w_k)^2}{(h+\delta_k \e^1\cdot \w_k)^2}}-1\right)= o(\delta_k^2)+\frac {\delta_k^2}{2}\frac{|\w_k|^2
-(\e^1\cdot \w_k)^2}{h+\delta_k \e^1\cdot \w_k}$$
and similarly
$$
\left| \h +\delta_k \g\right|-\e^1\cdot \left(\h+\delta_k \g \right)= o(\delta_k^2)+\frac{\delta_k^2}{2}\frac{|\g|^2-\left( \e^1\cdot \g\right)^2 }{h+\delta_k\e^1\cdot \g}.
$$

Letting $k\rightarrow \infty$ we may then 
drop the assumption that $\g$ is bounded.
In particular, for $\g$ such that $\g=\w$ in $B_1 \setminus B$, we arrive at the inequality
$$
\int_{B} |\nabla \w|^2 + \int_{B} \frac{|\w|^2-(\e^1\cdot \w)^2}{h}
\le
\int_{B} |\nabla \g|^2 + \int_{B} \frac{|\g|^2-(\e^1\cdot \g)^2}{h}
$$
for all $\g \in W^{1,2}\left( B_1;\mathbb{R}^m\right) $ coinciding with $\w$ on $\partial B$. 

Calculation of the first variation yields that
\begin{align*}
\Delta (\e^1\cdot \w)&=0 \quad  \text{in}  \quad B,\\
\Delta (\e^j\cdot \w)&= \frac{\e^j\cdot \w}{h} \quad\text{in}\quad B \quad \text{ for }\quad j>1.
\end{align*}
By Lemma \ref{eigen} as well as the homogeneity of $\w$ 
and the fact that $\w\equiv 0$ in $B_1^-$
we obtain 
$$\e^j\cdot \w(x) = d_j h(x) \textrm{ for each } j>1,$$
where $d_j$ is a constant real number.
\\{\bf Step 3: $w:= \e^1\cdot \w =0$ in $B_1$.}

As $w$ is harmonic in $B_1^+$, homogeneous of degree $2$ and satisfies $w=0$ in $B_1^-$ we obtain
(using for example odd reflection and the Liouville theorem)
that
$w(x) = \sum_{j=1}^{n-1} a_{nj} x_jx_n$ in $B_1^+$.
Remember that we have chosen $\h$ as the minimiser of
$\inf_{\h \in \Hf} 
\Vert \c_k-\h\Vert_{W^{1,2}(B_1;\R^m)}$.
It follows that for $\h_\nu := \e^1 \max(x\cdot \nu,0)^2/2$,
$$\frac{(\w_k,\h_\nu-\h)_{W^{1,2}(B_1;\R^m)}}{|\nu-\e^n|} \le
\frac{1}{2\delta_k} \frac{\Vert \h_\nu-\h\Vert^2_{W^{1,2}(B_1;\R^m)}}{|\nu-\e^n|}
\to 0\textrm{ as } \nu\to \e^n.$$
Therefore
\begin{align*}
&o(1) \ge \int_{B_1} \w_k \cdot \e^1 \Bigg[
\chi_{\{x_n>0 \} \cap \{x\cdot \nu>0\}} \frac{(x\cdot\nu)^2-(x\cdot \e^n)^2}{|\nu-\e^n|}
- 
\chi_{\{x_n>0 \} \cap \{x\cdot \nu\le 0\}} \frac{(x\cdot\e^n)^2}{|\nu-\e^n|}\\&
+
\chi_{\{x_n\le 0 \} \cap \{x\cdot \nu>0\}} \frac{(x\cdot\nu)^2}{|\nu-\e^n|}
\Bigg]\\
&+ 
\int_{B_1}  \Bigg[
\chi_{\{x_n>0 \} \cap \{x\cdot \nu>0\}} \frac{x\cdot(\nu+\e^n)(\nu-\e^n)+x\cdot(\nu-\e^n)(\nu+\e^n)}{|\nu-\e^n|}
\\&- 
\chi_{\{x_n>0 \} \cap \{x\cdot \nu\le 0\}} \frac{2x_n \e^n}{|\nu-\e^n|}
+
\chi_{\{x_n\le 0 \} \cap \{x\cdot \nu>0\}} \frac{2x\cdot \nu \nu}{|\nu-\e^n|}
\Bigg]\cdot \nabla \w_k \cdot \e^1.
\end{align*}
Setting $\xi := \lim_{\nu\to \e^n} \frac{\nu-\e^n}{|\nu-\e^n|}$, we see that for $\nu\to \e^n$
\begin{align*}&
\frac{(x\cdot\nu)^2-(x\cdot \e^n)^2}{|\nu-\e^n|} \to 2 x_n x\cdot \xi,
\frac{x\cdot(\nu+\e^n)(\nu-\e^n)}{|\nu-\e^n|}\to 2x_n \xi,
\frac{x\cdot(\nu-\e^n)(\nu+\e^n)}{|\nu-\e^n|}\to 2x\cdot \xi \e^n.
\end{align*}
On the other hand, on the set $(\{x_n>0 \} \cap \{x\cdot \nu\le 0\})
\cup(\{x_n\le 0 \} \cap \{x\cdot \nu>0\})$,
$|x \cdot \nu|=O(|\nu-\e^n|)$ and $|x \cdot \e^n|=O(|\nu-\e^n|)$ as $\nu\to \e^n$.
Passing first to the limit $\nu\to \e^n$ we conclude that
\begin{align*}
&o(1) \ge 2 \int_{B_1} [\w_k\cdot \e^1 x\cdot \xi \max(x_n,0)
+ (\max(x_n,0)\xi + \chi_{\{ x_n>0\}}x\cdot \xi \e^n)\cdot \nabla \w_k\cdot \e^1 ].
\end{align*}
Passing next to the limit $k\to\infty$, and taking into account that $\xi_n=0$ and that
$$\nabla w = \begin{pmatrix}a_{nj}x_n\\ \sum_{j=1}^{n-1}a_{nj}x_j\end{pmatrix},$$
we obtain that
\begin{align}\label{epi5}
&0\ge \sum_{j=1}^{n-1} a_{nj} \int_{B_1} [\max(x_n,0)^2 x\cdot \xi x_j + \max(x_n,0)^2 \xi_j
+ \chi_{\{x_n>0\}}x\cdot \xi x_j].
\end{align}
Since also $$\int_{B_1} x_j x_i=0  \qquad \textrm{ for } i\ne j,$$
we deduce from (\ref{epi5}) that
\begin{align}\label{epi6}
&0\ge \sum_{j=1}^{n-1} a_{nj} \xi_j \int_{B_1^+} (x_n^2x_j^2 + x_n^2 + x_j^2)
\qquad \textrm{ for every } \xi=(\xi_1,\dots,\xi_{n-1},0).
\end{align}
Thus $a_{nj}=0$ for $j=1,\dots,n-1$, that is $w\equiv 0$ in $B_1^+$.
\\
{\bf Step 4: $d_j=0$ for each $j\ge 2$.}

From Step 2-3 we know that $\w_k= \d h + \z_k$, where $\d\cdot \e^1=0$ and $\z_k\to 0$
weakly in $W^{1,2}(B_1;\R^m)$ as $k\to \infty$.
It follows that $\c_k = h(\e^1 +\delta_k \d) + \delta_k \z_k$.
By assumption,
\begin{equation}\label{epi20}
1 = \Vert \d h+  \z_k\Vert_{W^{1,2}(B_1;\R^m)}^2
= |\d|^2 \Vert h\Vert_{W^{1,2}(B_1;\R^m)}^2+2 (\d h,\z_k)_{W^{1,2}(B_1;\R^m)}
+ \Vert \z_k \Vert_{W^{1,2}(B_1;\R^m)}^2.
\end{equation}
Remember that we have chosen $\h$ as the minimiser of
$\inf_{\f \in \Hf} 
\Vert \c_k-\f \Vert_{W^{1,2}(B_1;\R^m)}$.
It follows that
for $\f := h(\e^1+\delta_k \d)/\sqrt{1+\delta_k^2 |\d|^2}\in \Hf$,
\begin{align*}
&\delta_k = \Vert \c_k - h \e^1\Vert_{W^{1,2}(B_1;\R^m)}
\le \Vert \c_k - \f\Vert_{W^{1,2}(B_1;\R^m)}
= \left\Vert h(\e^1 +\delta_k \d) + \delta_k \z_k - \frac{h(\e^1+\delta_k \d)}{\sqrt{1+\delta_k^2 |\d|^2}}
\right\Vert_{W^{1,2}(B_1;\R^m)}\\&
= \left\Vert \delta_k \z_k + \frac{h(\e^1+\delta_k \d)}{\sqrt{1+\delta_k^2 |\d|^2}}
(\sqrt{1+\delta_k^2 |\d|^2}-1)\right\Vert_{W^{1,2}(B_1;\R^m)}
\le \delta_k\Vert \z_k\Vert_{W^{1,2}(B_1;\R^m)} + C_3 \delta_k^2 |\d|^2.
\end{align*}
\begin{equation}\label{epi21}
  \textrm{Hence, }1\le \Vert \z_k\Vert_{W^{1,2}(B_1;\R^m)} + C_3 \delta_k|\d|^2.
\end{equation}
Combining (\ref{epi20}) and (\ref{epi21}), we obtain that
\begin{align*}
|\d|^2 \Vert h\Vert_{W^{1,2}(B_1;\R^m)}^2&+2 (\d h,\z_k)_{W^{1,2}(B_1;\R^m)}
+ \Vert \z_k \Vert_{W^{1,2}(B_1;\R^m)}^2\\
&\le \Vert \z_k\Vert_{W^{1,2}(B_1;\R^m)}^2+O(\delta_k).
\end{align*}
Letting $k\to\infty$, we conclude that
$|\d|^2 \Vert h\Vert_{W^{1,2}(B_1)}^2=0$ and that $|\d|=0$.
\\
{\bf Step 5: $\w_k \to \w$ strongly in $W^{1,2}(B_1;\R^m)$.}

Plug in $\v := (1-\zeta)\c_k + \zeta \h$ in (\ref{epi3}), where
$\zeta(x)=\min(2\max(1-|x|,0),1)$. Then
$$\frac{\v-\h}{\delta_k} = (1-\zeta) \w_k,$$ and we obtain that
\begin{align*}
&(1-\kappa_k) \bigg[\int_{B_1} |\nabla \w_k|^2  - 2 \int_{\partial B_1} 
|\w_k|^2 \haus
+ 2 \int_{B_1^-} \frac{|\c_k|}{\delta_k^2} + 2  \int_{B_1^+} \frac{|\c_k|- \e^1\cdot \c_k}{\delta_k^2}\bigg]
\\&<
\int_{B_1} |\nabla ((1-\zeta)\w_k)|^2  - 2 \int_{\partial B_1} 
|(1-\zeta)\w_k|^2 \haus
+ 2 \int_{B_1^-} \frac{(1-\zeta)|\c_k|}{\delta_k^2}\\& + 2  \int_{B_1^+} \frac{|(1-\zeta)\c_k+\zeta \h|}{\delta_k^2}
-2 \int_{B_1^+} \frac{(1-\zeta)\c_k\cdot \e^1+\zeta h}{\delta_k^2}.
\end{align*}
Using the definition of $\zeta$, it follows that
\begin{align*}
&\int_{B_{1/2}} |\nabla \w_k|^2 \le C_4 \kappa_k + \int_{B_1}\left(
 |\nabla \zeta|^2 |\w_k|^2 - 2 (1-\zeta) \nabla \zeta
\cdot \nabla \w_k \w_k\right).
\end{align*}
 The integral on the left-hand side equals by homogeneity of $\w_k$
$$2^{-n-2} \int_{B_1} |\nabla \w_k|^2,$$
so that
\begin{align*}
&\int_{B_1} |\nabla \w_k|^2 \le 2^{n+2}  
\left(C_4 \kappa_k + \int_{B_1}\left(
 |\nabla \zeta|^2 |\w_k|^2 - 2 (1-\zeta) \nabla \zeta
\cdot \nabla \w_k \cdot \w_k\right)\right)\to 0 \textrm{ as } k\to\infty.
\end{align*}
\\
{\em Altogether we obtain a contradiction from $\w\equiv 0$,
the strong convergence of $\w_k$ as well as the fact that $\Vert \w_k\Vert_{W^{1,2}(B_1;\R^m)}=1$.}
\qed


\section{Introduction to the problem and technical tools}\label{tools}

Let $D$ be a bounded Lipschitz domain in $\R^n$ and let
$\u=(u_1,\dots,u_m)$ be a minimiser of 
$$E(\u) := \int_D (|\nabla \u|^2 + 2 |\u|)$$
in the affine subspace $\{ \v\in W^{1,2}(D; \mathbb{R}^m): \v=\u_D
\textrm{ on } \partial D\}.$
Note that non-negativity, convexity and lower semicontinuity with respect to weak
convergence imply existence of a minimiser for each $\u_D\in W^{1,2}(D; \mathbb{R}^m)$.

In order to compute the first variation of the energy, we compute for $\vphi \in W^{1,2}_0\left( D ; \mathbb{R}^m\right) $
\begin{align}\label{varin}
0\le \epsilon \int_D 2 \nabla \u \cdot \nabla \vphi
+ \epsilon^2 \int_D |\nabla \vphi|^2 + 2 \int_D (|\u+\epsilon \vphi|-|\u|)\\
\le \epsilon \int_D 2 \nabla \u \cdot \nabla \vphi
+ \epsilon^2 \int_D |\nabla \vphi|^2 + 2|\epsilon| \int_D |\vphi|.\non
\end{align}
Dividing by $\epsilon$ and letting $\epsilon\to 0$, it follows that 
$$\left|\int_D  \nabla \u \cdot \nabla \vphi\right|\le  \Vert \vphi \Vert_{L^1(D;\R^m)},$$
so that $\Delta \u\in L^\infty(D; \mathbb{R}^m)$.
Applying standard $L^p$- and $C^\alpha$-theory, we obtain that
$\u\in W^{2,p}_{\rm loc}(D; \mathbb{R}^m)\cap C^{1,\alpha}_{\rm loc}(D; \mathbb{R}^m)$ for each $p\in [1,+\infty)$
and each $\alpha \in (0,1)$.
We see that $\Delta \u = 0$ a.e. in $\{ \u=0\}$.
Moreover, in the open set $\{ |\u|>\delta>0\}$,
passing to the limit in (\ref{varin})
yields $$\Delta \u = \frac{\u}{|\u|} \textrm{ in }  \{ |\u|>\delta>0\}.$$
Altogether we obtain
that $\u$ is a strong solution of the equation
$$\Delta \u = \frac{\u}{|\u|} \chi_{\{ |\u|>0\}}$$
in $D$.

Note that any other solution $\v\in W^{1,2}(D; \mathbb{R}^m)$ with the same boundary data $\u_D$ and
satisfying the weak equation
$$\int_D \left(\nabla \v \cdot \nabla \vphi + \vphi \frac{\v}{|\v|}\chi_{\{ |\v|>0\}}\right)=0
\textrm{ for every } \vphi\in W^{1,2}_0(D; \mathbb{R}^m)$$
must coincide with $\u$:
plugging in $\vphi := \v-\u$ yields
$$\int_D |\nabla (\u-\v)|^2
\le - \int_D \left(\frac{\u}{|\u|}\chi_{\{ |\u|>0\}}-\frac{\v}{|\v|}\chi_{\{ |\v|>0\}}\right)\cdot (\u-\v)\le 0.$$
Thus the weak solution is unique and equals the minimiser of the problem, so that
it is sufficient to consider minimisers. 

Note that in contrast to the classical (scalar) obstacle problem, it is an open problem
whether $\u\in W^{2,\infty}_{\rm loc}(D; \mathbb{R}^m)$.
\begin{remark}\label{stanest}
Using standard elliptic theory combined with the 
estimate $|\Delta \u|\le 1$ we obtain that
\begin{equation}
\sup_{B_{3/4}}|\u|+\sup_{B_{3/4}} |\nabla \u| \le C_1(n,m) \left(\Vert \u\Vert_{L^1(B_1;\R^m)}+1\right).\label{nablaest}
\end{equation}
\end{remark}
\begin{remark}
If a sequence of solutions of our system $\u_k$ converges weakly in $W^{1,2}(D;\R^m) $ 
to $\u$, then Rellich's theorem together
with the fact that $D^2 \u=0$ a.e. in $\{\u=0\}$, implies that $\u$ is a solution, too. 
\end{remark}
\begin{prop}[Non-Degeneracy]\label{ndeg}
Let $\u$ be a solution of  (\ref{system}) in $D$. If $x^0\in \overline{\{|\u|>0\}}$ and $B_r(x^0)\subset D$, then 
$$
\sup_{B_r(x^0)}|\u|\ge \frac{1}{2n}r^2.
$$
\end{prop}
\proof
It is sufficient to prove a uniform estimate for $x^0\in \{ |\u|>0\}$.
Let $U(x):= |\u (x)|$. Then
\begin{equation}\label{U}
\Delta U = 1+\frac{A}{U} \textrm{ in } \{|\u|>0\}
\textrm{, where } A=|\nabla \u|^2 - |\nabla U|^2 \ge 0.
\end{equation}
Assuming 
$\sup_{B_r(x^0)}|\u|\le \frac {1}{2n} r^2$, we obtain that the function
$$ v(x) := U(x)- U(x^0)-\frac{1}{2n}|x-x^0|^2$$
is subharmonic in the connected component of $ B_r(x^0)\cap \left\lbrace |\u|>0\right\rbrace $ containing
$x^0$, that $v<0$ on the boundary of that component and that $v(x^0)=0$, contradiction.\qed

\begin{prop}
Let $\u$ be a solution of  (\ref{system}) in $B_1(0)$ such that $\Vert \u-\h \Vert_{L^1(B_1;\R^m)}\le \epsilon<1$,
where $\h := \frac{\max(x_n,0)^2}{2}\e^1$. Then 
$$ B_{1/2}(0)\cap \supp \u \subset \left\lbrace x_n>-C\epsilon ^{\frac {1}{2n+2}}\right\rbrace $$ 
with a constant $C=C(n,m).$  
\end{prop}

\proof
Suppose that $B_{1/2}\cap \{ |\u|>0\}\ni x^0$ and that $x^0_n=-\rho<0$.
It follows that 
$$\Vert \u\Vert_{L^1(B_{\rho}(x^0);\R^m)}\le 
\Vert \u-\h\Vert_{L^1(B_1;\R^m)}\le \epsilon.$$
By the non-degeneracy property Proposition \ref{ndeg} we know that
$$|\u(y)|=\sup_{B_{\rho/2}(x^0)} |\u| \ge \frac{1}{8n} \rho^2$$ for some $y\in B_{\frac{\rho}{2}}(x^0)$.
From Remark \ref{stanest} we infer that  
$$\inf_{B_{\sigma\rho^2}(y)}|\u| \ge \frac{1}{8n} \rho^2 - 2C_1(n,m) \sigma \rho^2
\ge \frac{1}{16\> n}\rho^2,$$
provided that $\sigma$ has been chosen small enough, depending only on $n$ and $m$.
Combining our estimates, we obtain that
$$\epsilon\ge \Vert \u\Vert_{L^1(B_{\sigma\rho^2}(y);\R^m)}
\ge \left(\frac{1}{16\> n} \rho^2\right) |B_1| \left(\sigma\rho^2\right)^n,$$
a contradiction, if $\epsilon< C_2(n,m)\rho^{2n+2}$.  It follows that
$|\u(x)|=0$ for $x_n\le -C(n,m)\epsilon^{\frac{1}{2n+2}}$.
\qed


\section{Monotonicity Formula and Consequences}\label{mon-sec}

\begin{lem}\label{mon}
Let $\u$ be a solution of  (\ref{system}) in $B_{r_0}(x^0)$ and let 
$$
W(\u,x^0,r)=\frac{1}{r^{n+2}}\int_{B_r(x^0)}(|\nabla \u|^2+2|\u|)-
\frac{2}{r^{n+3}}\int_{\partial B_r(x^0)} |\u|^2\haus.
$$
For $0<r<r_0$,
$$
\frac{d W(\u,x^0,r)}{d r}=2\int_{\partial B_1(0)} r\left|\frac{d}{dr} \u_r\right|^2 \haus,
$$
where $\u_r(x)=\u(rx+x^0)/r^2$.
\end{lem}
\proof
The proof of this lemma follows by now standard arguments of G.S. Weiss (see \cite{cpde} and \cite{inv}).
A short proof consists in scaling
\begin{align*}&
\frac{d}{dr}W(\u,x^0,r)=
\frac{d}{dr}
\left(\int_{B_1(0)}(|\nabla \u_r|^2+2|\u_r|)-
2 \int_{\partial B_1(0)} |\u_r|^2\haus\right)\\&
= 2\int_{B_1(0)} (\nabla \u_r\cdot \frac{d}{dr}
\nabla \u_r + \frac{\u_r }{|\u_r|}\cdot \frac{d}{dr} \u_r)
- 2 \int_{\partial B_1(0)} 2\u_r \cdot \frac{d}{dr} \u_r\haus
\\&= \frac{2}{r}\Bigg(\int_{B_1(0)} (\nabla \u_r\cdot \nabla (x\cdot\nabla \u_r - 2\u_r)
 + \frac{ \u_r}{|\u_r|} \cdot (x\cdot\nabla \u_r - 2\u_r))\\
&\quad\quad
- 2 \int_{\partial B_1(0)} \u_r \cdot (x\cdot\nabla \u_r - 2\u_r)\haus\Bigg)\\&
= \frac{2}{r}\Bigg(\int_{B_1(0)} (-\Delta \u_r\cdot (x\cdot\nabla \u_r - 2\u_r)
 + \frac{\u_r }{|\u_r|}\cdot (x\cdot\nabla \u_r - 2\u_r))\\&\quad\quad
+\int_{\partial B_1(0)} (x\cdot\nabla \u_r-2\u_r) 
\cdot (x\cdot\nabla \u_r - 2\u_r)\haus\Bigg)\\&
= \frac{2}{r}\int_{\partial B_1(0)} |x\cdot\nabla \u_r - 2\u_r|^2\haus
= 2r\int_{\partial B_1(0)} \left| \frac{d \u_r}{dr}\right|^2\haus .
\end{align*}
This proves  the statement of the lemma.\qed

Note that for $x^0\in B_{1/2}$ and $r<1/2$,
\begin{equation}\label{energybound}
W(\u,x^0,r)\le C\big(\| \u \|_{W^{1,2}(B_1;\R^m)}+\|\u\|^2_{W^{1,2}(B_1;\R^m)}   \big).
\end{equation}
Moreover, we obtain the following properties:

\begin{lem}\label{properties}
1. The function $r\mapsto W(\u,x^0,r)$ has a right limit
$W(\u,x^0,0+)\in [-\infty,+\infty)$.

2. Let $0<r_k \to 0$ be a sequence such that the blow-up sequence
$$ \u_k(x) :=
\frac{\u(x^0+r_k x)}{r_k^2}$$
converges weakly in $W^{1,2}_{\rm loc}(\R^n;\R^m)$ to $\u_0$. 
Then $\u_0$ is a homogeneous function of degree $2$.
Moreover
 $$W(\u,x^0,0+) = \int_{B_1(0)} |\u_0| \ge 0,$$
and $W(\u,x^0,0+) =0$ implies that $u\equiv 0$ in $B_\delta(x^0)$
for some $\delta>0$.

3. The function $x \mapsto W(\u,x,0+)$ is upper-semicontinuous.
\end{lem}
\proof
1. follows directly from the monotonicity formula.

2. By the assumption of convergence $(\u_k)_{k\in \N}$ is
bounded in $W^{1,2}_{\rm loc}(\R^n;\R^m)$ and the limit $W(\u,x^0,0+)$ is finite.
From the monotonicity formula we obtain for all
$0<\rho<\sigma<+\infty$ that
$$\int_\rho^\sigma \frac{1}{r^{n+4}} \int_{\partial B_r(0)}|x\cdot\nabla \u_k(x) - 2 \u_k(x)|^2\haus\> dr
\to 0, k\to\infty,$$
proving the homogeneity of $\u_0$.

 We calculate, using the homogeneity of $\u_0$,
\begin{align*}&
W(\u,x^0,0+) = \int_{B_1(0)}(|\nabla \u_0|^2+2|\u_0|)-
2 \int_{\partial B_1(0)} |\u_0|^2\haus
\\&= \int_{B_1(0)}(-\u_0\cdot \Delta \u_0+2|\u_0|)+
\int_{\partial B_1(0)} (x\cdot\nabla \u_0\cdot\u_0 - 2|\u_0|^2)\haus
= \int_{B_1(0)} |\u_0| \ge 0.
\end{align*}
In the case $W(\u,x^0,0+) =0$ we obtain a contradiction to the non-degeneracy
Lemma \ref{ndeg} unless
$\u\equiv 0$ in some ball $B_\delta(x^0)$.

3. For $\epsilon>0$, $M<+\infty$ and $x\in D$ we obtain from the monotonicity
formula that
$$W(\u,x,0+) \le W(\u,x,\rho) \le \frac{\epsilon}{2}+W(\u,x^0,\rho)
\le \left\{\begin{array}{ll}
\epsilon+W(\u,x^0,0+),& W(\u,x^0,0+)>-\infty,\\
-M,& W(\u,x^0,0+)=-\infty,
\end{array}\right.$$
if we choose first $\rho$ and then $|x-x^0|$ small enough.
\qed


\section{A quadratic growth estimate}


\begin{thm}\label{quadratic}
Any solution $\u$ to the system (\ref{system}) in $B_1(0)$
satisfies 
$$
|\u (x)| \leq C\dist^2 (x,\Gamma_0(u))\textrm{ and }
|\nabla \u (x)| \leq C\dist(x,\Gamma_0(u)) \textrm{ for every }
x\in B_{1/2}(0),
$$
where the constant $C$ depends only on $n$ and
$$I(\u,0,1) :=
 \int_{B_1(0)}(|\nabla \u|^2+2|\u|).$$
\end{thm}
\proof

The statement of the theorem is equivalent to
$$
\sup_{B_r(x^0)}|\u|\le C_1r^2\textrm{ and } \sup_{B_r(x^0)}|\nabla\u|\le C_1r\textrm{ for every }
x^0\in \Gamma_0(\u)\cap B_{1/2}(0)\textrm{ and every } r\in (0,1/4),
$$
which in turn can be readily derived 
by standard elliptic theory from
\begin{equation}\label{intest}
\frac{1}{r^n}\int_{B_r(x^0)}|\u|\le C_2r^2\textrm{ for all } x^0 \textrm{ and } r \textrm{ as above.}
\end{equation}
Thus, our goal here is to show that (\ref{intest}) holds. To that end,
notice first that by the monotonicity formula,
$$W(\u,x^0,r)\le W(\u,x^0,1/2) \le 2^{n+2}I(\u,0,1)
\textrm{ for every } x^0\in B_{1/2}(0) \cap  \Gamma_0(\u) \textrm{ and } \ r\le 1/2.$$
Therefore
\begin{align*}
\frac{2}{r^{n+2}} \int_{B_r(x^0)} |\u|
&= W(\u,x^0,r)
- \frac{1}{r^{n+2}}\int_{B_r(x^0)}|\nabla \u|^2 + \frac{2}{r^{n+3}}\int_{\partial B_r(x^0)}|\u|^2\haus\\&
= W(\u,x^0,r) - \frac{1}{r^{n+2}}\int_{B_r(x^0)}|\nabla (\u-S_{x^0}\p)|^2 +  \frac{2}{r^{n+3}}\int_{\partial B_r(x^0)}|\u-S_{x^0}\p|^2\haus
\\&\le I(\u,0,1)
+ \frac{2}{r^{n+3}} \int_{\partial B_r(x^0)}|\u-S_{x^0}\p|^2\haus
\end{align*}
for each $\p=(p_1,\dots,p_m)\in \mathcal{H}$;
here the set $\mathcal{H}$ is the set
of all $\p=(p_1,\dots,p_m)$
such that each component $p_j$ is a 
homogeneous harmonic polynomial of second order, $S_{x^0}\f(x):=\f(x-x^0)$.

Let $x^0\in \Gamma_0$ and $\p_{x^0,r}$ be the 
minimiser of $\int_{\partial B_r(x^0)}|\u-S_{x^0}\p|^2\haus$
in $\mathcal{H}$.
It follows that
\begin{equation}\label{orth1}0=\int_{\partial B_r(x^0)}(\u-S_{x^0}\p_{x^0,r})\cdot S_{x^0}\q\haus
\textrm{ for every }\q\in \mathcal{H}.
\end{equation}
We maintain that 
there is a constant $C_1$ depending only on the dimension $n$ as well as
$I(\u,0,1)$ such that 
for each $x^0\in B_{1/2}(0)\cap \Gamma_0$ and $r\le 1/4$,
$$\frac{1}{r^{n+3}} \int_{\partial B_r(x^0)}|\u-S_{x^0}\p_{x^0,r}|^2\haus\le C_1.$$
Suppose towards a contradiction that there is a sequence of solutions
$\u_k$ (to equation (\ref{system}) in $B_1(0)$)  and a sequence of points $x^k
\in B_{1/2}(0)\cap \Gamma_0(\u_k)$ as well as $r_k\to 0$
such that $I (\u_k, 0, 1)$ are uniformly bounded, 
 $$M_k:=\frac{1}{r_k^{n+3}} \int_{\partial B_{r_k}(x^k)}|\u_k-S_{x^k}\p_{x^k,r_k}|^2\haus\to\infty \ .$$
For $\v_k(x) := \u_k(x^k+r_k x)/r_k^2$, 
and $\w_k(x) := \left( \v_k- \p_{x^k,r_k}\right) /M_k$, we have
 $\Vert \w_k\Vert_{L^2(\partial B_1(0);\R^m)}=1$
and
\begin{align*}
\int_{B_1(0)} |\nabla \w_k|^2 - 2 \int_{\partial B_1(0)} |\w_k|^2 
&= M_k^{-2} 
\left(
\int_{B_1(0)} 
|\nabla (\v_k- \p_{x^k,r_k})|^2 - 
2 \int_{\partial B_1(0)} |\v_k- \p_{x^k,r_k}|^2\haus
\right)\\
&= M_k^{-2} \left(\int_{B_1(0)} |\nabla \v_k|^2 - 
2 \int_{\partial B_1(0)} |\v_k|^2\haus\right)\\
&\le M_k^{-2} W(\u_k,x^k,r_k)\le M_k^{-2} C_2\to 0,\ k\to\infty. 
\end{align*}
It follows that $(\w_k)_{k\in \N}$ is bounded in $W^{1,2}(B_1;\R^m)$ such that
---passing to a subsequence if necessary--- $\w_k$ converges weakly
in $W^{1,2}(B_1;\R^m)$ to $\w_0\in W^{1,2}(B_1;\R^m)$.
By Rellich's theorem, $\w_k$ converges strongly
in $L^2(\partial B_1(0);\R^m)$, $\Vert \w_0\Vert_{L^2(\partial B_1(0);\R^m)}=1$
and (by (\ref{orth1}))
$0=\int_{\partial B_1(0)}\w_0\cdot \q\haus$
for every $\q\in \mathcal{H}$. Hence we obtain that
\begin{equation}\label{freq11}\int_{B_1(0)}|\nabla\w_0|^2\le 2\int_{\partial B_1(0)}|\w_0|^2\haus=2.\end{equation}
Moreover,
$$|\Delta \w_k|\le \frac{C_3}{M_k} \to 0,k\to\infty,$$
such that $\w_0$ is harmonic in $B_1(0)$ and
(by $C^{1,\alpha}$-estimates) $|\w_0|(0)=|\nabla \w_0|(0)=0$.
Now by \cite[Lemma 4.1]{interfaces}, 
each component $z_j$ of $\w_0$ must satisfy
 $$2 \int_{\partial B_1(0)} z_j^2\haus \le  \int_{B_1(0)} |\nabla z_j|^2.$$
Summing over $j$ we obtain
$$2 \int_{\partial B_1(0)} |\w_0|^2\haus \le  \int_{B_1(0)} |\nabla \w_0|^2,$$
implying by (\ref{freq11}) that
 $$2 \int_{\partial B_1(0)} |\w_0|^2\haus = \int_{B_1(0)} |\nabla \w_0|^2.$$
Thus 
$$2 \int_{\partial B_1(0)} z_j^2\haus =\int_{B_1(0)} |\nabla z_j|^2$$
for each $j$, implying by \cite[Lemma 4.1]{interfaces} that
$z_j$ is a homogeneous harmonic polynomial of second order.
But then
$0=\int_{\partial B_1(0)}\w_0\cdot \q\haus$
for every $\q\in \mathcal{H}$ implies that $\w_0=0$ on $\partial B_1(0)$,
contradicting $\Vert \w_0\Vert_{L^2(\partial B_1;\R^m)}=1$.\qed

The next section follows closely the procedure in \cite{inv} and \cite{jde}.


\section{An energy decay estimate and uniqueness of blow-up limits}\label{secdecay}


In this section we show that an epiperimetric inequality always
implies an energy decay estimate and uniqueness of blow-up limits.
More precisely:

\begin{thm}[Energy decay and uniqueness of blow-up limits] \label{decay}
Let $x^0\in D\cap\partial\{|\u|>0\}$, and 
suppose that the epiperimetric inequality
holds with $\kappa \in (0,1)$ for each 
$$
\c_r(x) := {\vert x\vert}^2
\u_r({x\over {\vert x \vert}}) = {{\vert x \vert}^2 \over {r^2}}
\u(x^0+{r\over {\vert x \vert}} x)
$$
and for all $r\le r_0 < 1$. Finally  let $\u_0$ denote an arbitrary blow-up limit of $\u$ at $x^0 $.
Then 
$$
\left\vert W(\u,x^0,r)-W(\u,x^0,0+)\right\vert
\; \le \; \left\vert W(\u,x^0,r_0)-W(\u,x^0,0+)\right\vert
\> \left({r\over {r_0}}\right)^{(n+2)\kappa\over {1-\kappa}}
$$
for $r \in (0,r_0) $, and  there exists a constant $C$ depending
only on $n$ and $\kappa$ such that
\[ \int_{\partial B_1(0)} \left\vert {\u(x^0+rx)\over {r^2}} 
\> - \> \u_0(x)\right\vert \> d\mathcal{H}^{n-1} \; \le \; 
C \>  \left\vert W(\u,x^0,r_0)-W(\u,x^0,0+)\right\vert^{1\over 2}
\> \left({r\over {r_0}}\right)^{(n+2)\kappa\over {2(1-\kappa)}}\]
for $r \in (0,{r_0\over 2})\> ,$ and $\u_0$ is the unique
blow-up limit of $\u$ at $x^0\> .$
\end{thm}

\proof
We define
\[ e(r) := r^{-n-2} \int_{B_r(x^0)} ({\vert \nabla \u \vert}^2 \> + \>
2|\u|)\; - \; 2 \> r^{-n-3} \> \int_{\partial B_r(x^0)}
|\u|^2 \> d\mathcal{H}^{n-1}\; - \; W(\u,x^0,0+).\]
Up to a constant $e(r)$ is the function of the monotonicity identity,
so that we have already computed $e'(r)$. Here however, we need
a different formula for $e'(r)$:
\[ e'(r) \; = \; \Bigg[ \> -{n+2\over r}
e(r) \; - \; {n+2\over r} W(\u,x^0,0+)\; +
\; {2\over r} \> 
r^{-n-3} \int_{\partial B_r(x^0)} |\u|^2 
\> d\mathcal{H}^{n-1}\]
\[ - \; 2 \> r^{-n-3} \> \int_{\partial B_r(x^0)}
2 \nu \cdot\nabla\cdot\u \u \> d\mathcal{H}^{n-1} \; 
- \; {2(n-1)\over r} \> r^{-n-3} \int_{\partial B_r(x^0)} |\u|^2 
\> d\mathcal{H}^{n-1} \] \[ + \> r^{-n-2}\> \int_{\partial B_r(x^0)}
({\vert \nabla \u \vert}^2 \> + \>
2|\u|) \> d\mathcal{H}^{n-1}\Bigg]
\; = \;
r^{-1}\Bigg[ \int_{\partial B_1(0)} ({\vert \nabla \u_r\vert}^2
\] \[ + \> |\u_r| \> - \> 4 \nu\cdot\nabla \u_r\cdot\u_r\> + \> 4 |\u_r|^2
\> + \> 4 |\u_r|^2
\> - \> 2(n+2)|\u_r|^2) \> d\mathcal{H}^{n-1}
\; - \; (n+2) W(\u,x^0,0+)\Bigg] \] \[- \; {n+2\over r} e(r)
\; \ge \; r^{-1}\bigg[ \int_{\partial B_1(0)} ({\vert \nabla_\theta
\u_r \vert}^2 \> + \> |\u_r| \> + \> 4|\u_r|^2\> - \> 2(n+2)|\u_r|^2
) \> d\mathcal{H}^{n-1} 
\] \[- \; (n+2) W(\u,x^0,0+)\bigg] \; - \; {n+2\over r} e(r)
\; = \; r^{-1}\bigg[ \int_{\partial B_1(0)} ({\vert \nabla_\theta
\c_r \vert}^2 \> + \> |\c_r| \> + \> |\nu \cdot\nabla \c_r |^2 )
\> d\mathcal{H}^{n-1} \] \[ - \; (n+2)2 \> \int_{\partial B_1(0)}
{\c_r}^2 \> d\mathcal{H}^{n-1} \; - \; (n+2) W(\u,x^0,0+)\bigg]
\; - \; {n+2\over r} e(r)
\] \[= \; {n+2\over r} \left[ M(\c_r) -  W(\u,x^0,0+)- e(r)\right] \; .\]
Here we employ the minimality of $\u$ as well as the assumption
that the epiperimetric inequality $M(\v) \le (1-\kappa) M(\c_r)
\> + \> \kappa \> W(\u,x^0,0+)$ holds for some $\v \in W^{1,2}(B_1;\R^m)$
with $\c_r$-boundary values and we obtain 
for $r\in (0,r_0)$ the estimate
\[ e'(r) \ge {n+2\over r} {1\over {1-\kappa}} (M(\u_r)-W(\u,x^0,0+))
\; - \; {n+2\over r} e(r) \] \[= \; {n+2\over r} \left(
{1\over {1-\kappa}} - 1\right) e(r) \; = \; {(n+2)\kappa\over
{1-\kappa}} {1\over r} \> e(r)\; .\]
By the monotonicity formula Lemma \ref{mon}, $e(r)\ge 0\> ,$ and we conclude
in the non-trivial case $e>0$ in $(r_1,r_0)$ that
\[ (\log(e(s)))' \; \ge \; {(n+2)\kappa\over
{1-\kappa}} {1\over s} \; \hbox{ for } s\in (r_1,r_0)\; .\]
Integrating from $r$ to $r_0$ we obtain that
\[ \log\left({e(r_0)\over {e(r)}}\right) \; \ge \;
{(n+2)\kappa\over
{1-\kappa}} \log\left({r_0\over r}\right)\; \hbox{ and } \; 
{e(r_0)\over {e(r)}} \; \ge \; \left({r_0\over r}\right)^{(n+2)\kappa\over
{1-\kappa}} \; \hbox{ for } r \in (r_1,r_0)\]
and that $e(r) \le e(r_0) \> \left({r \over {r_0}}\right)^{(n+2)\kappa\over
{1-\kappa}}$ for $r \in (0,r_0)$ which proves our first
statement.\\
Using once more the monotonicity formula (Lemma \ref{mon}) we get
for $0<\rho<\sigma\le r_0$ an estimate of the form
\[ 
\int_{\partial B_1(0) } \int_\rho^\sigma 
 \left\vert \frac{d \u}{dr} \right\vert dr 
 \> d\mathcal{H}^{n-1}  \le \;
\int_\rho^\sigma \int_{\partial B_1(0)} r^{-2} 
\left\vert {\nu \cdot\nabla \u(x^0+rx)} \> - \> 2 {\u(x^0+rx)\over {r}}
\right\vert \> d\mathcal{H}^{n-1} \> dr
\] \[= \; \int_\rho^\sigma r^{-1-n} \int_{\partial B_r(x^0)} 
\left\vert \nu \cdot\nabla \u \> - \> 2 {\u\over r}\right\vert 
\> d\mathcal{H}^{n-1}
\> dr\] \[\le \; \sqrt{n \> \omega_n}  \int_\rho^\sigma
r^{-1-n} r^{n-1\over 2} r^{n+2\over 2}\left(
r^{-n-2} \int_{\partial B_r(x^0)} 
{\left\vert \nu \cdot\nabla \u \> - \> 2 {\u\over r}\right\vert}^2 \> d\mathcal{H}^{n-1}
\right)^{1\over 2} 
\> dr \] \[= \; \sqrt{n \> \omega_n\over 2}\int_\rho^\sigma
r^{-{1\over 2}}\> \sqrt{e'(r)} \> dr \;\le \; \sqrt{n \> \omega_n\over 2}
(\log(\sigma)-\log(\rho))^{1\over 2}
\> (e(\sigma)-e(\rho))^{1\over 2}\; . \]
Considering now $0<2\rho<2r\le r_0$ and intervals
$[2^{-k-1},2^{-k}) \ni \rho$ and $[2^{-\ell-1},2^{-\ell})\ni r$ the already
proved part of the theorem yields that
\[ \int_{\partial B_1(0) }\left\vert
{\u(x^0+rx)\over {r^2}} \> - \> {\u(x^0+\rho x)\over {\rho^2}}
\right\vert \> d\mathcal{H}^{n-1} \leq 
\sum_{i=\ell}^k
\int_{\partial B_1(0) }\int_{ 2^{-i-1}}^{2^{-i}} 
\left\vert 
\frac{d \u}{dr}\right\vert dr 
\> d\mathcal{H}^{n-1}
\] 
 \[\le \> 
C_1(n) \sum_{i=\ell}^k (\log(2^{-i})-\log(2^{-i-1}))^{1\over 2}
(e(2^{-i})-e(2^{-i-1}))^{1\over 2} \> = \>
C_2(n) \sum_{i=\ell}^k (e(2^{-i})-e(2^{-i-1}))^{1\over 2} \] \[\le \;
C_3(n,\kappa) \left\vert W(\u,x^0,r_0)-W(\u,x^0,0+)\right\vert^{1\over 2}
\sum_{j=\ell}^{+\infty} \left(r_0\> 2^j\right)^{-(n+2)\kappa\over
{2(1-\kappa)}} \] \[\le \; 
C_4(n,\kappa) \left\vert W(\u,x^0,r_0)-W(\u,x^0,0+)\right\vert^{1\over 2} 
\> \frac{c^\ell}{1-c}
\> {r_0}^{-(n+2)\kappa\over
{2(1-\kappa)}}\] where $c=2^{-(n+2)\kappa\over
{2(1-\kappa)}} \in (0,1)\> .$ Thus
\[ \int_{\partial B_1(0) }\left\vert
{\u(x^0+rx)\over {r^2}} \> - \> {\u(x^0+\rho x)\over {\rho^2}}
\right\vert \> d\mathcal{H}^{n-1} \] \[\le \; 
C_5(n,\kappa) \left\vert W(\u,x^0,r_0)-W(\u,x^0,0+)\right\vert^{1\over 2}
\left({r\over {r_0}}\right)^{(n+2)\kappa\over
{2(1-\kappa)}}\; ,\] and letting ${\u(x^0+\rho_j x)\over {{\rho_j}^2}}\to \u_0$
as a certain sequence $\rho_j\to 0$ finishes our proof.\quad\qed


\section{Homogeneous solutions}


In this section we consider homogeneous solutions $\u \in W^{1,2}(B_1;\R^m) $, meaning 
that
$$
\u (\lambda x)=\lambda^2 \u(x) \quad \textrm{ for all } \lambda >0 \textrm{ and }x\in B_1(0).
$$
Obviously $\u$ may be extended to a homogeneous solution on $\mathbb{R}^n$. 

Moreover, if $\dist_{L^1(B_1;\R^m)}(\u,\Hf)\le 1$ then due to Remark 1 we have
\begin{equation}\label{standest2}
\sup_{B_1} |\u| \le C(n,m)\textrm{ and } \sup_{B_1} |\nabla \u| \le C(n,m).
\end{equation}

\begin{prop} \label{homog}
If $B_1\cap \supp \u\subset \{ x_n >- \delta(n,\sup_{B_1(0)}|\u|)\}$, then
$\u \in \Hf$.
\end{prop}
\proof
Observe first that each component $u_i$ is a solution of 
$$\mathcal{L}_iu_i := -\Delta' u_i + u_i/|\u| = 2n u_i$$
in every connected component $\Omega'$ of $\partial B_1\cap \{ |\u|>0\}$,
where $\Delta'$ is the Laplace-Beltrami operator on the unit sphere in
$\R^n$. In Lemma \ref{eigen} of the Appendix we prove that 
$u_i = a_if_{\Omega'}$ for a real number $a_i$ and a function $f_{\Omega'}$ 
depending only on $\Omega'$ which
is positive on $\Omega'$ and vanishes on the boundary of $\Omega'$.
It follows that for each connected component $\Omega$ of $B_1\cap \left\lbrace |\u|>0\right\rbrace $
there exists a unit vector $\a=(a_1\dots,a_m)$ such that
$\u(x)=\a|\u(x)|$ and $\Delta |\u|=1$ in $\Omega$.

Now, if $|\nabla \u|=0$ on $\partial \Omega$, then 
we may extend $\u$ by $0$ outside $\Omega$,
that is $|\u|$ can be extended to 
a $2$-homogeneous non-negative solution of the classical
obstacle problem in $\R^n$. These solutions have been completely classified (see \cite{caff-1981}, cf. also \cite{psu}), and
$\supp \u\subset \left\lbrace  x_n >- \delta |x|\right\rbrace $ (where $\delta=\delta(n,m,\sup_{B_1(0)}|\u|)$) would in this case imply 
that up to rotation, $|\u|=h$, and $\u=\a h$, where $h$ is  a half-space solution
for scalar problem. 

If, on the other hand, there is a point $x^0\in \partial\Omega\cap
\{|\nabla \u|\ne 0\}$, then the fact that $\u$ is continuously differentiable,
implies that $\a$ equals the vector
of the adjoining connected component of $\{ |\u|>0\}$ up to the sign.
In this case we obtain, taking the maximal union of all such connected components,
that each $u_i$ is a $2$-homogeneous solution of the scalar two-phase obstacle problem
$$\Delta v = c(\chi_{\{v>0\}} - \chi_{\{v<0\}}) \textrm{ in }\R^n$$ with $c>0$,
satisfying $v=0$ in $\{ x_n\le -\delta\}$.
However, according to \cite[Theorem 4.3]{monats}, no such solution exists.
\qed
\begin{lem}\label{isol}
The half-plane solutions are (in the $L^{1}(B_1(0);\R^m)$-topology)
isolated within the class of homogeneous
solutions of degree $2$.
\end{lem}
\proof

Let $\Vert \u-\h\Vert_{L^1(B_1;\R^m)}\le \epsilon$, where 
rotating in $\R^n$ and in $\R^m$ if necessary we may assume that 
$$\h(x)= \frac{\max(x_n,0)^2}{2} \e^1.$$ From (\ref{standest2}) as well as Propositions 2 and 3
we infer that $\u \in \Hf$ if $\epsilon$ has been chosen small enough, depending only on 
$n$ and $m$.
\qed

We defined earlier the constant
$ \alpha_n=2M(\h)$ where $\h \in \Hf$. Now we are going to estimate the value
of $M(\u)$ for an arbitrary homogeneous solution $\u$  of degree 2.

\label{sopen}
\begin{prop}
\begin{equation}\label{glob1}
\alpha_n = \frac{\mathcal{H}^{n-1}(\partial B_1)}{2n(n+2)}.
\end{equation}
Let $\u$ be a homogeneous solution of degree $2$. Then
\begin{equation}\label{glob2}
M(\u) \ge \alpha_n \frac{\mathcal{H}^{n-1}(\partial B_1\cap
\{ |\u|>0\})}{\mathcal{H}^{n-1}(\partial B_1)}.
\end{equation}
In particular,
\begin{align}\label{glob3}
&M(\u) \ge \alpha_n \textrm{ if } |\u|>0 \textrm{ a.e.}
\end{align}
\end{prop}
\proof
Let $U:= |\u|$, and recall (\ref{U}):
$$\Delta U = 1+\frac{A}{U} \textrm{ in } \{|\u|>0\}
\textrm{, where } A=|\nabla \u|^2 - |\nabla U|^2 \ge 0.$$
It follows that ---using the homogeneity of $\u$---
\begin{align*}
&\int_{B_1\cap \{|\u|>0\}} \left(1+\frac{A}{U}\right)
= \int_{\partial (B_1\cap \{|\u|>0\})} \nabla U \cdot \nu \haus
\\&=  2\int_{\partial B_1\cap \{|\u|>0\}} U \haus- 2\int_{B_1\cap \partial\{ |\u|=0\}\cap \{|\nabla \u|> 0\}}
|\nabla U|\haus
\\&= 2(n+2) \int_{B_1} U - 2\int_{B_1\cap \partial\{ \u=0\}\cap \{|\nabla \u|> 0\}}
|\nabla U|\haus.\end{align*}
On the other hand, using once more the homogeneity of $\u$,
\begin{align}\label{sobl1}
M(\u) = \int_{B_1} (|\nabla \u|^2 + 2 |\u|) -  2 \int_{\partial B_1} |\u|^2 \haus=
\int_{B_1} |\u| = \int_{B_1} U.
\end{align}
In order to verify (\ref{glob1}),
observe that for $\e\in \partial B_1 \subset \R^m$ and $\h(x)=\e \max(x_n,0)^2/2$,
\begin{align*}
\frac{\alpha_n}{2}=M(\h) = \frac{1}{2} \int_{B_1} \max(x_n,0)^2
= \frac{1}{4} \int_{B_1} x_n^2
=\frac{1}{4n} \int_{B_1} |x|^2
= \frac{\mathcal{H}^{n-1}(\partial B_1)}{4n(n+2)}.
\end{align*}
Using the above estimates we conclude that
\begin{align*}
M(\u) \ge \frac{1}{2(n+2)} 
|B_1\cap
\{ \u>0\}|
= \frac{\mathcal{H}^{n-1}(\partial B_1\cap
\{ \u>0\})}{2n(n+2)} = \alpha_n \frac{\mathcal{H}^{n-1}(\partial B_1\cap
\{ \u>0\})}{\mathcal{H}^{n-1}(\partial B_1)}.\end{align*} \qed

\begin{cor}\label{min}
Let $\u$ be a homogeneous solution of degree $2$. Then
\begin{equation}\label{glob4}
M(\u)\ge \alpha_n \max \left( \frac{1}{2},\frac{\mathcal{H}^{n-1}(\partial B_1\cap
\{ |\u|>0\})}{\mathcal{H}^{n-1}(\partial B_1)}\right) \ge \alpha_n/2,
\end{equation}
and $M(\u)=\alpha_n/2$ implies that $\u\in \Hf$.
Moreover, $\alpha_n/2<\bar\alpha_n := \inf \{ M(\v): \v$ is a
homogeneous solution of degree $2$, but $\v\not\in \Hf\}$.
\end{cor}
\proof
If $\mathcal{H}^{n-1}(\partial B_1\cap
\{ |\u|=0\})=0$, then \eqref{glob4} follows from (\ref{glob3}).
Otherwise $\{ |\u|=0\}$ contains by the non-degeneracy property Lemma \ref{ndeg}
an open ball $B_\rho (y)$, and we may choose it in such a way  that there is a point
$z\in \partial B_\rho (y)\cap
\partial \{ |\u|>0\}$. Let $\u_0$ be a blow up of $\u$ at $z$. Since $\supp \  \u_0$ is
contained in a half-space it follows from Proposition~\ref{homog} that
$\u_0 \in \Hf$.
From (\ref{glob1}) we obtain therefore that
\begin{equation}\label{homeq}
\frac{\alpha_n}{2}=  W(\u,z,0^+)\le W(\u,z,+\infty) = W(\u,0,+\infty)
= M(\u).
\end{equation}

Now we have to prove that $M (u)= \alpha_n /2$ implies $\u \in \Hf$.
Consider a ball $B_\rho (y)$ and a point $z$ as above,  that is $y=z+\rho e$ with a unit vector $e$.
It follows from homogeneity of $\u$ that $e$ is orthogonal to $z$. We consider two cases.

\noindent
{\bf Case a)} If $z=0$ then, again due to homogeneity of $\u$, we have $|\u(x)|=0$
in a half-space  $(x\cdot e)>0$. Hence $\u\in \Hf$ by Proposition 3.

\noindent
{\bf Case b)} If  $|z|>0$ then, since $W(\u,z,r)$ does not depend on $r$ (by (\ref{homeq})), we conclude 
that $u$ is homogeneous with homogeneity  center $z$. More exactly, we have  
 $\u(z+kx)=k^2\u(z+x)$ for any $k>0, x\in \mathbb{R}^n$. Since also $\u(z+kx)=
k^2\u(z/k+x)$ we obtain  $u(z+x)=u(z/k+x)$
for any $k>0, x\in \mathbb{R}^n$. It means that $\u$ is constant in direction
of vector $z$.  In particular, $|\u|=0$ in the ball 
$B_\rho(\rho e)$ touching the origin and we are again at the case a). 

Last, we have to prove $\bar{\alpha}_n>\alpha_n/2$. 
If it is not true then there is a sequence of homogeneous global solutions $\left\lbrace \u_k\right\rbrace $ such that
$$
M(\u_k)\searrow \frac{\alpha_n}{2} \quad \textrm{as} \quad k\rightarrow \infty.
$$
In particular it implies by (\ref{sobl1}) uniform boundedness of $\u_k$ in $L^1(B_1(0))$ and
therefore, by  (\ref{nablaest}) and by elliptic theory uniform boundedness of solutions $\u_k$  
in $W^{2,q}_{\rm loc}\left( \mathbb{R}^n\right) $ for any $q<\infty$. Then there exists a limit $\tilde{\u}$, by subsequence, such that
$\tilde{\u}$ is a homogeneous solution, $\tilde{\u} \not \in \Hf$ (by Lemma \ref{isol}) and $M(\tilde{\u})=\alpha_n/2$. From the first part of the proof we infer that $\tilde{\u} \in \Hf$, and a 
contradiction arrises.
\quad\qed


\begin{definition}\label{regul-bdry}
A point $x$ is a regular free boundary point for $u$ if: 
$$ x\in  \Gamma_0 (\u)    \qquad  \hbox{ and }\qquad 
\lim_{r\to 0} W(\u,x,r) = {\alpha_n \over 2}   .
$$
We denote by $\mathcal{R}_u$ the set of all regular free boundary points of $u$ in $B_1$.
\end{definition}

\begin{cor}\label{op}
The set of regular free boundary points $\mathcal{R}_u$
is open relative to $\Gamma_0(\u).$
\end{cor}
\proof
This is an immediate consequence of  Corollary
\ref{min} and the upper semicontinuity Lemma \ref{properties}.\quad\qed





\section{Regularity}\label{sregular}

In this last section we prove that the set of regular free boundary
points $\mathcal{R}_u$ is locally in $D$ a $C^{1,\beta}$-surface and we
derive a {\em macroscopic criterion for regularity:}
suppose that $W(\u,x,r)$ drops for some (not necessarily small)
$r$ below the critical value $\bar \alpha_n\> :$ then
$\partial\{|\u|>0\}$ must be a $C^{1,\beta}$-surface in an open 
neighborhood of $x\> .$


\begin{thm}  \label{tdiff}
Let $C_h$ be a compact set of points $x^0\in \Gamma_0(\u)$ with the following property: at least one 
blow-up limit $\u_0$ of $\u$ at $x^0$ is a half-plane solution,
say $\u_0(x) = {1\over 2} \e \max(x\cdot \nu(x^0),0)^2$ for some
$\nu(x^0)\in \partial B_1(0)\subset \R^n$ and $\e(x^0) \in \partial B_1(0)\subset \R^m$.
Then there exist $r_0>0$ and $C<\infty$ such that 
\[ \int_{\partial B_1(0)} \left\vert {\u(x^0+rx)\over {r^2}}
\> - \> {1\over 2} \e(x^0)\max(x\cdot \nu(x^0),0)^2\right\vert
\> d\mathcal{H}^{n-1} \; \le \; C \> r^{(n+2)\kappa\over {2(1-\kappa)}}\] \[
\hbox{ for every } x^0\in C_h \hbox{ and every } r \le r_0\; .\]
\end{thm}
\proof
In view of Theorem \ref{decay} and Theorem \ref{epi},
it is sufficient to show that
$$
\dist_{W^{1,2}(\partial B_1(0);\R^m)} \left( {\u(x^0+r\cdot)\over {r^2}},\Hf\right) 
+ \dist_{L^\infty(\partial B_1;\R^m)} \left( {\u(x^0+r\cdot)\over {r^2}},\Hf\right)  \le \delta
$$
for every $x^0\in C_h$ and $r\le 2r_0.$

To this end, we first prove that for $\epsilon>0$ there exists
$\tilde\delta>0$ such that the estimate
\begin{equation}\label{mod}
W(\u,x,r)-{\alpha_n\over 2} \le \epsilon
\end{equation}
holds for every $x \in C_h$ and $r\le \tilde\delta\> :$
suppose that this is not true; then we find sequences $r_k\to 0\> , \>
C_h \ni x^k\to x^0$ as $k\to \infty$
such that $W(\u,x^k,r_k)-{\alpha_n\over 2}\; > \; \epsilon\> .$ 
By the monotonicity of $W $, for $\rho \geq r_k$ we can deduce
\[ \epsilon < W(\u,x^k,r_k)-{\alpha_n\over 2}  \leq
  W(\u,x^k,\rho)  - W(\u,x^0,\rho) +  W(\u,x^0,\rho)
-  W(\u,x^0,0+) < \epsilon ,
\] 
provided we choose first $\rho$ small and then $k$ large.\\
Considering now {\em any} sequences ${\u(x^k+\rho_k\cdot) \over
{{\rho_k}^2}}$ such that $x^k\in C_h$ and
$\rho_k\to 0$ as $k \to \infty\> ,$
the monotonicity formula Lemma \ref{mon}, (\ref{mod}) and the fact
that ${\u(x^k+\rho_k\cdot) \over
{{\rho_k}^2}}$ is bounded in $W^{2,p}_{\rm loc}(\R^n;\R^m)$ imply
that all limits $\u_h$ of ${\u(x^k+\rho_k\cdot) \over
{{\rho_k}^2}}$ must be homogeneous solutions
of degree $2$ satisfying $M(\u_h)={\alpha_n \over 2}$.\\
Assume now towards a contradiction that 
$\dist_{W^{1,2}(\partial B_1(0);\R^m)}({\u(x^k+\rho_k \cdot)\over {
{\rho_k}^2}},\Hf)\ge \delta > 0$ 
as the sequence $k\to \infty$. Since ${\u(x^k+\tau_\ell \cdot) \over
{{\tau_\ell}^2}}$ converges by assumption for a certain sequence 
$\tau_\ell\to 0$ to some $h\in \Hf$ as $\ell\to \infty\> ,$
we find by a continuity argument for each $\theta \in (0,1)$ a sequence $
\tilde\rho_k\to 0$ such that $ 
\dist_{W^{1,2}(\partial B_1(0);\R^m)}({\u(x^k+\tilde\rho_k \cdot)\over {
{\tilde\rho_k}^2}},\Hf)= \theta\delta\> .$ Passing to a limit
with respect to a subsequence we obtain a homogeneous solution 
$\u_h$ of degree $2$ satisfying 
$\dist_{W^{1,2}(\partial B_1(0);\R^m)}(\u_h,\Hf)= \theta\delta$ which
contradicts the isolation property Lemma \ref{isol}.\quad\qed

\begin{thm}[regularity] \label{regular}
The set of regular free boundary points 
$\mathcal{R}_u$ is locally in $D$ a $C^{1,\beta}$-surface; here
$\beta = {(n+2)\kappa\over {2(1-\kappa)}} \left(
1\> + \> {(n+2)\kappa\over {2(1-\kappa)}}\right)^{-1}\; .$
\end{thm}
\proof
Let us consider a point $x^0 \in \mathcal{R}_u\> .$ By Corollary
\ref{op} and Theorem \ref{tdiff} there exists $\delta_0>0$
such that $B_{2\delta_0}(x^0)\subset D\> ,\>
B_{2\delta_0}(x^0)\cap \partial\{|\u|>0\} = B_{2\delta_0}(x^0)\cap \mathcal{R}_u$
and 
\begin{equation}\label{est1}
\int_{\partial B_1(0)} \left\vert {\u(x^1+rx)\over {r^2}}
\> - \> {1\over 2} \e(x^1)\max(x\cdot \nu(x^1),0)^2\right\vert
\> d\mathcal{H}^{n-1} \; \le \; C \> r^{(n+2)\kappa\over {2(1-\kappa)}}
\end{equation}\[
\; \hbox{ for every } x^1\in \partial\{|\u|>0\}\cap\overline{B_{\delta_0}(x^0)} 
\hbox{ and for every } r \le \min(\delta_0,r_0)\; .\]
We now observe that $x^1\mapsto \nu(x^1)$ and $x^1\mapsto \e(x^1)$ are H\"older-continuous
with exponent $\beta$ on $\partial\{|\u|>0\} \cap
\overline{B_{\delta_1}(x^0)}$ for some $\delta_1\in (0,\delta_0)\> :$ 
\[ {1\over 2} \int_{\partial B_1(0)}\left\vert\e(x^1)\max(x\cdot \nu(x^1),0)^2
\> - \> \e(x^2)\max(x\cdot \nu(x^2),0)^2\right\vert \> d\mathcal{H}^{n-1} 
\; \le \; 2C\> r^{(n+2)\kappa\over {2(1-\kappa)}} \] \[ + \;
\int_{\partial B_1(0)} \int_0^1 \left\vert {\nabla \u(x^1+rx+t(x^2-x^1))
\over {r^2}} \right\vert \> \vert x^1-x^2\vert\>
dt \> d\mathcal{H}^{n-1} \; \le \;
2C\> r^{(n+2)\kappa\over {2(1-\kappa)}} \] \[+ \; C_1 \> {\max(r,\vert
x^1-x^2\vert)\vert x^1-x^2\vert\over {r^2}}
\; \le \; (2C+C_1) {\vert x^2-x^1\vert}^{\gamma
{(n+2)\kappa\over {2(1-\kappa)}}}\]
if we choose $\gamma := \left(1\> + \> {(n+2)\kappa\over {2(1-\kappa)}}\right)^
{-1}$ and $r:= {\vert x^2-x^1\vert}^{\gamma} \le \min(\delta_0,r_0)\> ,$
and the left-hand side 
\begin{align*}&{1\over 2} \int_{\partial B_1(0)}\left\vert\e(x^1)\max(x\cdot \nu(x^1),0)^2
\> - \> \e(x^2)\max(x\cdot \nu(x^2),0)^2\right\vert \> d\mathcal{H}^{n-1}
\\& \ge \; c(n) \> (\vert \nu(x^1)-\nu(x^2)\vert + \vert \e(x^1)-\e(x^2)\vert)\end{align*}
as can be easily shown by an indirect argument.\\
Next, (\ref{est1}) as well as the regularity and non-degeneracy
of $\u$ imply that for $\epsilon>0$ there exists $\delta_2\in (0,\delta_1)$
such that 
\begin{equation}\label{cone}
\begin{array}{l}
\u(y) = 0 \hbox{ for } x^1 \in \partial\{|\u|>0\} \cap
\overline{B_{\delta_1}(x^0)} \\
\hbox{ and } y \in 
\overline{B_{\delta_2}(x^1) } \hbox{ satisfying }
(y-x^1)\cdot \nu(x^1) < -\epsilon\vert y-x^1\vert \hbox{ and }\\
|\u(y)| > 0 \hbox{ for } x^1 \in \partial\{|\u|>0\} \cap
\overline{B_{\delta_1}(x^0)} \\
\hbox{ and } y \in 
\overline{B_{\delta_2}(x^1)} \hbox{ satisfying }
(y-x^1)\cdot \nu(x^1) > \epsilon \vert y-x^1 \vert  \    : 
\end{array}
\end{equation}
assuming that (\ref{cone}) does not hold, we obtain a sequence 
$\partial\{|\u|>0\} \cap \overline{B_{\delta_1}(x^0)}\ni x^m \to \bar x$
and a sequence $y^m-x^m\to 0$ as $m\to \infty$ such that 
\begin{equation}
\label{ann}
\begin{array}{c}
\hbox{either } |\u(y^m)|>0 \hbox{ and } 
(y^m-x^m)\cdot \nu(x^m) < -\epsilon\vert y^m-x^m\vert\\
\hbox{or } \u(y^m)=0 \hbox{ and } 
(y^m-x^m)\cdot \nu(x^m) > \epsilon\vert y^m-x^m\vert\; .\end{array}
\end{equation}
On the other hand we know from (\ref{est1}) as well as from the
regularity and non-degeneracy of the solution $\u$, that the sequence
$\u_j(x) := {\u(x^j+\vert y^j-x^j\vert \> x)\over {{\vert y^j
- x^j \vert}^2}}$ converges in $C^{1,\alpha}_{\rm loc}(\R^n;\R^m)$
to ${1\over 2}\e(\bar x)\max(x\cdot \nu(\bar x),0)^2$ as $j \to \infty$ and that
$\u_j=0$ on each compact subset $C$ of $\{x\cdot \nu(\bar x)<0\}$
provided that $j \ge j(C)\> .$ This, however, 
contradicts (\ref{ann}) for large $j\> .$\\
Last, we use (\ref{cone}) in order to show that $\partial\{|\u|>0\}$
is for some $\delta_3 \in (0,\delta_2)$
in $\overline{B_{\delta_3}(x^0)}$ the graph of a differentiable
function: applying two rotations we may assume that $\nu(x^0)=\e^n$ and $\e(x^0)=\e^1$. 
Choosing now $\delta_2$ with respect to $\epsilon={1\over 2}$
and defining functions $g^+,g^-:B'_{\delta_2\over 2}(0)\to [-\infty,
\infty]\> , \> g^+(x') := \sup\{ x^n\> : \> x^0\> + \> (x',x^n)
\in \partial\{|\u|>0\} \}$ and $g^-(x') := \inf\{ x^n\> : \> x^0\> + \> (x',x^n)
\in \partial\{|\u|>0\} \}$ we infer from (\ref{cone}) as well as
from the continuity of $\nu(x)$ immediately that 
$\{ x^n\> : \> x^0\> + \> (x',x^n)
\in \partial\{|\u|>0\} \}$ is non-empty and that for sufficiently
small $\delta_3$ the functions $g^+$ and $g^-$ are Lipschitz-continuous
and satisfy $g^+=g^-$ on $\overline{B'_{\delta_3}(0)}\> .$
Applying (\ref{cone}) once more with respect to arbitrary $\epsilon$
we see that $g^+$ is Fr\'echet-differentiable in 
$\overline{B'_{\delta_3}(0)}\> ,$ which finishes our proof in view
of the already derived H\"older-continuity of the normal $\nu(x)\> .$
\quad\qed

\begin{cor}[Macroscopic criterion for regularity]
\label{crit}
There exists $\bar \alpha_n > {\alpha_n \over 2}$ such that 
$B_{2r}(x^0)\subset D\> , \> x^0 \in D\cap\partial\{|\u|>0\}$
and $W(\u,x^0,r) < \bar \alpha_n$ imply that $\partial\{|\u|>0\}$
is in an open neighborhood of $x^0$ a $C^{1,\beta}$-surface.
\end{cor}
\proof
The first statement follows immediately from Lemma \ref{isol},
Corollary \ref{min} and
Theorem \ref{regular}. \qed

  

\section{Appendix}


\begin{lem}\label{eigen}
Let $\Delta'$ be the Laplace-Beltrami operator on the unit sphere in
$\R^n$, let the domain $ \Omega'\subset \partial B_1(0)\subset \R^n$,
let $\mathcal{L} := -\Delta' +q$ where $q\in C^0(\Omega')$ such that
$q\ge q_0>0$ in $\Omega'$,
and let $\lambda_k(\mathcal{L},\Omega')$ denote the $k$-th eigenvalue
with respect to the eigenvalue problem
\begin{align*}
&\mathcal{L} v = \lambda v \textrm{ in } \Omega'\\
&v=0 \textrm{ on }\partial \Omega';
\end{align*}
here $\partial \Omega'$ denotes the boundary of $\Omega'$ relative to $\partial B_1$.
\\
1. If $\tilde\Omega'\subset \Omega'$ then $\lambda_k(\mathcal{L},\tilde\Omega')\ge \lambda_k(\mathcal{L},\Omega')$ for every $k\in \N$. For $k=1$ the inequality is strict.
\\
2. $\lambda_k(\mathcal{L},\Omega')\ge q_0 + \lambda_k(-\Delta',\Omega')$ for every $k\in \N$;
in case $q\not\equiv q_0$ the inequality becomes a strict inequality.
\\
3. $q=1/h$ and $\Omega'\subset \partial B_1\cap \{ x^n >0\}$ and $v\in W^{1,2}(\partial B_1)$ being an eigenfunction with respect to $\Omega'$ and $\lambda=2n$ imply $v= a h$ for some real number $a\ne 0$. Here $h(x)=\frac{1}{2} \max\left(x_n,0\right)^2.$
\\
4. $\Omega' \subset \partial B_1(0) \cap \{ x^n > -\delta(n,q_0)\}$ and $v\in W^{1,2}(\partial B_1(0))$ being an eigenfunction of $\mathcal{L}$ with respect to $\lambda=2n$ and a domain $\Omega'$
imply $v= a f_{\Omega'}$ for a real number
$a\ne 0$ and a function $f_{\Omega'}$ which
is positive on $\Omega'$ and depends only on $\Omega'$. 
\end{lem}

\proof
1. It suffices to remark that $v\in W^{1,2}_0(\tilde\Omega')$ implies
 $v\in W^{1,2}_0(\Omega')$, after extending $v$ by zero outside $\tilde \Omega'$.

2. Let $\mathcal{M}_{k-1}$be a subspace of $W^{1,2}_0(\Omega')$ of codimension $k-1$,
$$ \mu(\mathcal{L},\Omega',\mathcal{M}_{k-1})\>=\>\inf_{v\in \mathcal{M}_{k-1},
\Vert v\Vert_{L^2(\Omega')}=1} \int_{\Omega'} (|\nabla' v|^2+qv^2).$$
Due to the Courant minimax principle we have 
$$\lambda_k(\mathcal{L},\Omega')\>=\>\sup\>\mu (\mathcal{L},\Omega',\mathcal{M}_{k-1})$$
where $\sup$ is taken over the set of all possible $\mathcal{M}_{k-1}$.
Since 
$$\mu(\mathcal{L},\Omega',\mathcal{M}_{k-1}) \ge \inf_{v\in \mathcal{M}_{k-1},\Vert v\Vert_{L^2(\Omega')}=1} 
\int_{\Omega'} |\nabla' v|^2+ \inf_{v\in \mathcal{M}_{k-1};\Vert v\Vert_{L^2(\Omega')}=1} 
\int_{\Omega'}qv^2
\ge q_0 + \mu(-\Delta',\Omega',\mathcal{M}_{k-1}),$$
we may take $\mathcal{M}_{k-1} := \left\lbrace v\in W^{1,2}_0(\Omega')\>:\>\int_\Omega' vw_i=0\right\rbrace $, where $w_i$ is an eigenfunction
with respect to the $i$th eigenvalue of $-\Delta'$ on $\Omega'$. For such $\mathcal{M}_{k-1}$
we obtain 
$ \mu(-\Delta',\Omega',\mathcal{M}_{k-1}) = \lambda_k(-\Delta',\Omega')$, and 2. is proved.

3. 
In this case the eigenvalue problem on the sphere becomes
\begin{align*}
&-\Delta' v + \frac{2}{\cos^2 \theta} v = 2n v \textrm{ in } \Omega'\subset \{ x^1>0\},
\\
&v=0 \textrm{ on }\partial \Omega'.
\end{align*}
Since $2n = \lambda_2(-\Delta',\partial B_1(0)\cap \{ x^n>0\})$,
we obtain from 1. and 2. that
$\lambda_2(\mathcal{L},\Omega') > 2+2n$.
But then $v$ must be an eigenfunction with respect to the {\em first}
eigenvalue $\lambda_1(\mathcal{L},\Omega')$, and $\lambda_1(\mathcal{L},\Omega')=2n$.
Observe now that $h$ is an eigenfunction with respect to the first
eigenvalue $\lambda_1(\mathcal{L},\partial B_1(0)\cap \{ x^n>0\})$, so that 
$$\lambda_1(\mathcal{L},\partial B_1(0)\cap \{ x^n>0\})=2n=\lambda_1(\mathcal{L},\Omega').$$
This is only possible if
$\Omega' = \partial B_1(0)\cap \{ x^n>0\}$ and 
$v = a h$ for some $a \ne 0$. 

4. Observe that the second eigenvalue of $-\Delta'$  in a half-sphere is $2n$, therefore
due to 1. and continuity of $\lambda_2$ with respect to the size of a spherical cap,
 $\lambda_2(-\Delta',\Omega') \ge  2n -\omega(n,\delta)$  where $\omega(n,\delta) \to 0$ as $\delta \to 0$. From 2. it follows that $\lambda_2(\mathcal{L},\Omega') \ge  2n -\omega(n,\delta)+q_0>2n$ if $\delta=\delta(n,q_0)$ is small. Thus, $v$ must be
 an eigenfunction with respect to the {\em first}
eigenvalue $\lambda_1(\mathcal{L},\Omega')$, and $\lambda_1(\mathcal{L},\Omega')=2n$.
Again, the eigenspace is a one-dimensional space such that
$v= a f_{\Omega'}$ for some real number
$a\ne 0$ and a function $f_{\Omega'}$ depending only on $\Omega'$;
We also know that first eigenfunctions do not change sign in the connected set $\Omega'$.
\qed

\bibliographystyle{plain.bst}
\bibliography{asuw-2013-07-30.bib}

\end{document}